\documentclass[12pt,a4paper,reqno]{amsart}
\usepackage{amsmath,amsthm}
\usepackage{amsfonts}
\usepackage{amssymb}
\usepackage{tikz-cd}
\allowdisplaybreaks

\newcommand{\be}{\begin{equation}}
\newcommand{\ef}{\end{equation}}
\chardef\bslash=`\\ 

\newtheorem*{thm*}{Theorem}

\theoremstyle{definition}

\newtheorem*{remark*}{Remarks}
\newtheorem*{defn*}{Definition}

\theoremstyle{remark}

\newcommand{\wt}{\widetilde}
\newcommand{\wh}{\widehat}
\newcommand{\fc}{\frac}

\newcommand{\iy}{\infty}
 \renewcommand{\sectionmark}[1]{}

\newcommand{\sh} {subharmonic}

\newcommand{\Te} {Teichm\"{u}ller}

\newcommand{\const}{\operatorname{const}}

 \usepackage{amsfonts}

\newcommand{\D}{\mathbb{D}}
\newcommand{\om}{\omega}
\newcommand{\z}{\zeta}
\newcommand{\ov}{\overline}
\newcommand{\vp}{\varphi}
\newcommand{\hC}{\wh{\mathbb{C}}}
\newcommand{\C}{\mathbb{C}}

\newcommand{\B}{\mathbf{B}}
\newcommand{\T}{\mathbf{T}}
\newcommand{\Belt}{\operatorname{Belt}}

\newcommand{\Om} {\Omega}
\newcommand{\vk} {\varkappa}
\newcommand{\x} {\mathbf x}
\renewcommand{\a} {\alpha}

\newcommand{\ld}{\lambda}

\begin{document}

\title{The Grunsky norm of univalent functions and abelian holomorphic differentials}
\author{Samuel L. Krushkal}

\begin{abstract} We establish that the Grunsky norm of any normalized univalent function on the disk
is completely determined by the squares of holomorphic abelian differentials (in contrast to the Teichm\"{u}ller
norm, which relates to all integrable holomorphic quadratic differentials).

This result has important interesting applications. In particular, it provides an explicit representation of Fredholm eigevalues of all quasiconformal curves.
\end{abstract}

\date{\today\hskip4mm({GruAbel(2).tex})}

\maketitle

\bigskip

{\small {\textbf {2020 Mathematics Subject Classification:} Primary: 30C55, 30C62, 30F45; Secondary 30F60, 32F45, 32G15, 46G20}

\medskip

}

\bigskip

\markboth{Samuel L. Krushkal}{The Grunsky norm of univalent functions and abelian holomorphic differentials}
\pagestyle{headings}

\bigskip\bigskip
\centerline{\bf 1. PRELIMINARIES}

\bigskip\noindent
{\bf 1.1. The Grunsky norm and its dominante}.
The classical Grunsky theorem of 1939 implies the necessary and sufficient conditions for univalence of holomorphic
functions in a finitely connected domain on the Riemann sphere $\hC = \C \cup \{\iy\}$ in terms of
an infinite system of the coefficient inequalities. In particular, for the canonical disk
$\D^* = \{z \in \hC: \ |z| > 1\}$ this theorem yields that a holomorphic function $f(z) = z + \const + O(z^{-1})$
in a neighborhood of $z = \iy$ can be extended to a univalent holomorphic function on $\D^*$ if and only if
the Taylor coefficients $\a_{mn}$ of the function
 \be\label{1}
\log \fc{f(z) - f(\z)}{z - \z} = - \sum\limits_{m,n=1}^\iy \a_{mn} z^{-m} \z^{-n}, \quad (z, \z) \in (\D^*)^2,
\end{equation}
called the {\bf Grunsky coefficients} of $f(z)$, satisfy the inequality
 \be\label{2}
\Big\vert \sum\limits_{m,n=1}^\iy \ \sqrt{m n} \ \a_{mn} x_m x_n \Big\vert \le 1,
\end{equation}
for any sequence $\mathbf x = (x_n)$ from the unit sphere $S(l^2)$ of the Hilbert space $l^2$ with norm $\|\x\| = (\sum\limits_1^\iy |x_n|^2)^{1/2}$; here the principal branch of the logarithmic function is chosen (cf. \cite{Gr}).
The quantity
$$
\vk(f) = \sup \Big\{\Big\vert \sum\limits_{m,n=1}^\iy \ \sqrt{mn} \ \a_{mn} x_m x_n \Big\vert: \
\mathbf x = (x_n) \in S(l^2) \Big\} \le 1
$$
is called the {\bf Grunsky norm} of $f$.

The univalent functions  $f(z) = z + b_0 + b_1 z^{-1} + \dots$ in $\D^*$ admitting quasiconformal extensions
across the unit circle $\mathbb S^1 = \partial \D^*$ onto the disk $\D = \{|z| < 1\}$ form the class $\Sigma_Q$.
To have their uniqueness for a given Beltrami coefficient $\mu(z) = \partial_{\ov z} f/\partial_z f$ in $\D$, compactness in the topology of locally uniform convergence on $\C$, etc., we add the third normalization condition
$$
f(0) = 0.
$$

All such  $f \in \Sigma_Q$ are zero free in $\D^*$, hence their inversions $F_f(z) = 1/f(1/z) = z + a_2 z^2 + \dots$
are holomorphic and univalent in the disk $\D$ with $F_f(\iy) = \iy$. The functions $f$ and $F_f$ have the same
Grunsky coefficients, and $\vk(F_f) = \vk(f)$.

Note also that the norm $\vk(f)$ is defined for all $f \in \Sigma_Q$ and does not depend on the additional normalization
at $0$.

For the functions with $k$-quasiconformal extensions ($k < 1$), we have instead of (2) a stronger bound
\be\label{3}
\Big\vert \sum\limits_{m,n=1}^\iy \ \sqrt{mn} \ \a_{mn} x_m x_n \Big\vert \le k \quad \text{for any} \ \ \x = (x_n) \in S(l^2),
\end{equation}
established first in \cite{Ku1} (see also \cite{Kr6}).

Note that the Grunsky matrix operator
$\mathcal G(f) = (\sqrt{m n} \ \a_{mn}(f))_{m,n=1}^\iy$
acts as a linear operator $l^2 \to l^2$ contracting the norms of elements $\x \in l^2$; the norm of this
operator equals $\vk(f)$.

The method of Grunsky inequalities was generalized in several directions, with the corresponding generalization of Grunsky norm.
In this paper continuing \cite{Kr7}, we shall deal only with the canonical case of disk $\D^*$.

\bigskip\noindent
{\bf 1.2. Some known results}.
The Grunsky norm $\vk(f)$ is dominated by the {\bf Teichm\"{u}ller norm} $k(f)$, which is equal to
the infimum of dilatations $k(w^\mu)= \|\mu\|_\iy$ of quasiconformal extensions of $f$ to $\hC$.
Here $w^\mu$ denotes a homeomorphic solution to the Beltrami equation $\partial_{\ov z} w = \mu \partial_z w$ on $\C$
extending $f$; accordingly, $\mu$ is called the {\bf Beltrami coefficient} (or complex dilatation) of $w$.

For most functions $f$, we have the strong inequality $\vk(f) < k(f)$ (moreover, the functions satisfying this inequality form a dense subset of $\Sigma_Q$), while the functions with the equal norms play a crucial role in many applications.

On the other hand, the important result of Pommerenke and Zhuravlev states  that if {\it a function $f \in \Sigma$ satisfies the inequality $\vk(f) < k$ with some constant $k < 1$, then $f$ has a quasiconformal extension to $\hC$ with a dilatation} $k_1 = k_1(k) \ge k$ \cite{Po}, \cite{Zh}; \cite{KK}, pp.82-84.

Each coefficient $\a_{m n}(f)$ in (4) is represented as a polynomial of a finite number of the initial coefficients
$b_1, b_2, \dots, b_s$ of $f$; hence it depends holomorphically on Beltrami  coefficients of quasiconformal extensions
of $f$ as well as on the {\bf Schwarzian derivatives}
$$
S_f(z) = \Bigl(\fc{f^{\prime\prime}(z)}{f^\prime(z)}\Bigr)^\prime - \fc{1}{2} \Bigl(\fc{f^{\prime\prime}(z)}{f^\prime(z)}\Bigr)^2, \quad z \in \D^*.
$$
These derivatives range over a bounded domain in the complex Banach space $\B(\D^*)$ of hyperbolically bounded holomorphic  functions $\vp \in \D^*$ with norm
$$
\|\vp\|_\B = \sup_{\D^*} (|z|^2 - 1)^2 |\vp(z)|.
$$
This domain models the {\bf universal Teichm\"{u}ller space} $\T$ (the space of complex structures on the disk)
in holomorphic Bers' embedding of $\T$. The complex geometry of this space provides a background of many results
related to this paper.

\bigskip
The Beltrami coefficients of quasiconformal extensions $w^\mu$ of functions $f(z) \in \Sigma_Q$ range over
the unit ball
$$
\Belt(\D)_1 = \{\mu \in L_\iy(\C): \ \mu(z)|\D^* = 0, \ \ \|\mu\|_\iy < 1\},
$$
and the well-known criterion for extremality (the Hamilton-Krushkal-Reich-Strebel theorem) implies that a Betrami coefficient
$\mu_0 \in \Belt(\D)_1$ is extremal if an only if
$$
\|\mu_0\|_\iy = \sup_{\|\psi\|_{A_1(\D)}=1} \ |\langle \mu, \psi\rangle_\D|,
$$
where
$$
\langle \mu, \psi\rangle_\D = \iint\limits_\D \mu(z) \psi(z) dx dy \quad (z = x + iy),
$$
and
$$
A_1(\D) = \{\psi \in L_1(D): \ \psi \ \ \text{holomorphic \ in} \ \ \D\}.
$$
The same condition is necessary and sufficient for the infinitesimal extremality of $\mu_0$ at the origin of the space $\T$
in the direction $t \phi_\T(\mu_0)$, where $\phi_\T$ is the defining (factorizing) holomorphic projection $\Belt(\D)_1 \to \T$; see, e.g., \cite{EKK}, \cite{GL}, \cite{Kr1}.

The Grunsky norm is connected with the subset of squares of the abelian holomorphic differentials $\om(z) dz$ on $\D$:
$$
A_1^2(\D) = \{\psi = \om^2 \in A_1(\D): \ \om \ \ \text{holomorphic \ in} \ \ \D\}
$$
formed by integrable squares of  squares of the abelian holomorphic differentials $\om(z) dz$.
Due to \cite{Kr2} (and to the more general result in  \cite{Kr6}), the equality $\vk(f^\mu) = k(f^\mu)$ is valid if and
only if
$$
\|\mu\|_\iy = \sup_{\psi \in A_1^2(\D),\|\psi\|_{A_1} =1} \ |\langle \mu, \psi\rangle_\D|;
$$
in addition, if $\mu$ is of Teichm\"{u}ller form, i.e., $\mu = k |\psi|/\psi$ with $\psi \in A_1(\D)$, then
necessarily $\psi \in A_1^2$.

Note also that, due to \cite{Kr2}, the elements of $A_1^2(\D)$ are represented in the form
$$
\psi(z) = \om(z)^2 = \fc{1}{\pi} \sum\limits_{m+n=2}^\iy \ \sqrt{mn} \ x_m x_n z^{m+n-2},
$$
with $\|\x\|_{l^2} = \|\om\|_{L_2}$; here $\x = (x_n)$.

\bigskip\noindent
{\bf 1.3. A lower bound for Grunsky norm}. The following essential estimate established in \cite{Kr7} is one of
the crucial steps in the proof of our main theorem.

\bigskip\noindent
{\bf Lemma 1}. {\it The Grunsky norm $\vk(f)$ of every function
$$
f(z) = z + b_0 + b_m z^{-m} + b_{m+1} z^{-(m+1)} + \dots \in \Sigma_Q, \quad m \ge 1
$$
satisfies the inequality
 \be\label{4}
\vk(f) \ge \a_\D(f) := \sup_{\psi \in A_1^2(\D),\|\psi\|_{A_1} =1}  \ |\langle\mu_0, \psi\rangle_\D|
= \sup_{\om \in A_2(\D), \|\om\|_2 =1} \ |\langle\mu_0, \om^2 \rangle_\D|,
\end{equation}
where $\mu_0$ is an extremal Beltrami coefficient among quasiconformal extensions $f^\mu$ of $f$ onto $\D$.  }

\bigskip
For functions $f$ with equal Grunsky and Teichm\"{u}ller norms, we have in (4) the equality, and vice versa
(see, \cite{Kr2}, \cite{Kr6}, \cite{Ku3}), but until now there was unknown whether there exist the functions
$f \in \Sigma_Q$ with $\vk(f) < k(f)$ and $\vk(f) = \a_\D(f)$.

The proof of this lower bound for $\vk(f)$ is based on geometric extensions of the Ahlfors Schwarz lemma to metrics
of integrable generalized negative curvatures given in \cite{Kr4}, \cite{Ro}.

We shall also use the infinitesimal version of the inequality (4). It will be presented later.

\bigskip\noindent
{\bf 1.4. Conformal metrics of negative generalized Gaussian curvature}.
We shall use some known results on subharmonic conformal metrics $ds = \ld(t) |dt|$ on the disk $\D$ (with $\ld(t) \ge 0$)
of negative generalized Gaussian curvature.

Recall that the {\bf generalized Gaussian curvature} $\kappa_\ld$ of an upper semicontinuous Finsler metric
$ds = \ld|dt|$ in a domain $\Om \subset \C$ is defined by
 \be\label{5}
\kappa_\ld (t) = - \fc{\Delta \log \ld(t)}{\ld(t)^2},
\end{equation}
where $\Delta$ is the {\it generalized Laplacian}
$$
\Delta \ld(t) = 4 \liminf\limits_{r \to 0} \frac{1}{r^2} \Big\{ \frac{1}{2 \pi} \int_0^{2\pi} \ld(t + re^{i \theta})
d \theta - \ld(t) \Big\}
$$
(provided that $- \iy \le \ld(t) < \iy$). Similar to $C^2$ functions, for which $\Delta$ coincides with the usual Laplacian,
one obtains that $\ld$ is \sh \ on $\Om$ if and only if
$\Delta \ld(t) \ge 0$; hence, at the points $t_0$ of local maximuma of $\ld$ with $\ld(t_0) > - \iy$, we have
$\Delta \ld(t_0) \le 0$.

The sectional {\it holomorphic curvature} of a Finsler metric on a complex Banach manifold $X$ is defined in a similar
way as the supremum of the curvatures (5) over appropriate collections of
holomorphic maps from the disk into $X$ for a given tangent direction in the image.

As is well-known \cite{AP}, \cite{Kr3}, the holomorphic curvature of the Kobayashi-Teichm\"{u}ller metric
$\mathcal K_\T(x, v)$ of
universal Teichm\"{u}ller space $\T$ equals $- 4$ at all points $(x, v)$ of the tangent bundle $\mathcal T(\T)$ over $\T$.
Instead, the holomorphic curvature of metric $\lambda_\vk$ generated on $\D$ by the Grunsky Finsler structure
satisfies the inequality
$\Delta \log \ld \ge 4 \ld^2$,
where $\Delta$ is again the generalized Laplacian (see \cite{Kr5}).
So, $\ld_\vk(t) \le -4$ on any holomorphic disk in the space $\T$.

This is a special case of metrics whose curvature satisfies the inequality $\Delta \log \ld \ge K \ld^2$
with constant $K > 0$. All such metrics are subharmonic.

We shall also apply the upper semicontinuous metrics $\ld$ satisfying the inequality
 \be\label{6}
\Delta \log \ld \ge K \ld, \quad K = \const > 0
\end{equation}
(then $u = \log \ld$ can be negative). For such metrics, we have the following {\bf Minda's maximum principle} given by

\bigskip\noindent
{\bf Lemma 2}. {\it If a function $u : \ D \to [- \iy, + \iy)$ is upper semicontinuous in a domain
$D \subset \C$ and its generalized Laplacian satisfies the inequality $\Delta u(z) \ge K u(z)$ with some positive
constant $K$ at any point $z \in D$, where $u(z) > - \iy$, and if
$$
\limsup\limits_{z \to \z} u(z) \le 0 \ \ \text{for all} \ \z \in \partial D,
$$
then either $u(z) < 0$ for all $z \in D$ or $u(z) \equiv 0$ on $D$.  }

\bigskip
The proof of this lemma related to the Ahlfors-Schwarz lemma is given in \cite{Mi}; its variations  see in \cite{Di}.

\bigskip\bigskip
\centerline{\bf 2. MAIN THEOREM AND ITS APPLICATIONS}

\bigskip\noindent
{\bf 2.1. Main theorem}.
The following theorem implies a complete answer on relation between the Grunsky norm and abelian holomorphic differentials.

\bigskip\noindent
{\bf Theorem 1}. {\it For any univalent function we have the equality}
 \be\label{7}
\vk(f) =  \a_\D(f).
\end{equation}

\bigskip
It suffices to consider the functions $f \in \Sigma_Q$, because for functions not admitting quasiconformal extension,
both sides of (7) are equal to $1$.

\bigskip\bigskip
{\bf 2.2. Applications}.
We mention two important applications of Theorem 1.

\bigskip\noindent
{\it 1. Plurisubharmonicity}. Letting $\a_\D(S_{f^\mu}) = \a_\D(f^\mu)$, one can regard this quantity as a function on
the universal Teichm\"{u}ller space $\T$ (as well as on the ball $\Belt(\D)_1$. The following corollary to Theorem 1 solves
the question posed in \cite{Kr4}.

\bigskip\noindent
{\bf Corollary 1}. {\it The function $\a_\D(S_{f^\mu})$ is plurisubharmonic on the space $\T$ and on the ball} $\Belt(\D)_1$.

\bigskip
This follows from the equality (7), because its left-hand side is plurisubharmonic on $\T$.
In fact, the indicated plurisubharmonicity is a consequence of the relation (10) and of holomorphy of functions (9)
defined below. Another proof is given in \cite{Kr6}.

\bigskip\noindent
{\it 2. Application to Fredholm eigenvalues of Jordan curves}.
The Fredholm eigenvalues $\rho_n$ of an oriented smooth closed Jordan curve $L  \subset \hC$ are the
eigenvalues of its double-layer potential, or equivalently, of the integral equation
$$
u(z) +  \frac{\rho}{\pi} \int\limits_L \ u(\zeta) \frac{\partial}{\partial
n_\zeta} \log \frac{1}{|\zeta - z|} ds_\zeta = h(z),
$$
where $n_\zeta$ denotes the outer normal and $ds_\zeta$ is the length element at $\zeta \in L$.
These values have crucial applications in solving many problems in various fields of mathematics.

It is important often to know the least positive eigenvalue $\rho_L = \rho_1$. This value is naturally connected with conformal and quasiconformal maps related to $L$ and can be defined for any oriented closed Jordan curve $L$ by
$$
\fc{1}{\rho_L} = \sup \ \fc{|\mathcal D_G (u) - \mathcal D_{G^*} (u)|} {\mathcal D_G (u) + \mathcal D_{G^*} (u)},
$$
where $G$ and $G^*$ are, respectively, the interior and exterior of $L; \ \mathcal D$ denotes the Dirichlet integral,
and the supremum is taken over all functions $u$ continuous on $\hC$ and harmonic on $G \cup G^*$. In particular,
$\rho_L = \iy$ only for the circle.

Note that both sides of the last equality remain invariant under the action of the Moebius group $PSL(2, \hC)$.

The indicated first eigenvalue is intrinsically connected with the Grunsky coefficients of the exterior conformal mapping
function  $f_L^*: \ \D^* \to D^*$;  this is qualitatively expressed by the remarkable K\"{u}hnau-Schiffer theorem, which
states that the value $\rho_L$ is reciprocal to the Grunsky norm $\vk(f_L^*)$ (see \cite{Ku2}, \cite{Sc}).

Together with this result, Theorem 1 implies explicitly the first Fredholm eigenvalue $\rho_L$ of every quasiconformal curve $L \subset \hC$. Namely, we have

\bigskip\noindent
{\bf Corollary 2}. {\it  For any quasicircle $L \subset \hC$,
  \be\label{8}
\fc{1}{\rho_L} = \sup_{\psi \in A_1^2(D),\|\psi\|_{A_1} =1} |\langle \mu_0, \psi\rangle_\D|,
\end{equation}
where $\mu_0$ is either of extremal Beltrami coefficients of the appropriately normalized exterior conformal mapping function $f_L^*$ on which the Teichm\"{u}ller norm of $f_L^*$ is attained.  }

Moreover, letting in (8) $\mu_0 = t |\psi|/\psi$ with $|t| < 1$ and $\psi \in A_1$, one obtains a dense subset of
a possible eigenvalues $\rho_L$.

\bigskip
Conversely, given a quasicircle $L \subset \hC$, take a quasiconformal extension $f^\mu$ of the exterior conformal map
$f_L^*$. Its Beltrami coefficient $\mu$ determines a linear functional
$\mu(\psi) = \langle \mu, \psi\rangle_\D$ on the linear span $\wt A$ of the subset $\{\psi \in A_1^2(\D): \ \|\psi\|_{A_1} = 1\}$
of the unit sphere in $A_1(\D)$. Its Hahn-Banach extension onto $L_1(\D)$ provides an extremal Beltrami coefficient $\mu_0$ in $\D$ for $f_L^*$ satisfying (8).

\bigskip\noindent
{\it 3. Teichm\"{u}ller maps with $\vk(f) = k(f)$}. We also mention the following useful property of Teichm\"{u}ller
extremal maps.

\bigskip \noindent
{\bf Corollary 3}. {\it If the equivalence class of $f$ is a Strebel point and $\vk(f) = k(f)$, then $\mu_0$ is
necessarily of the form
$\mu_0(z) = \|\mu_0\|_\iy |\psi_0(z)|/\psi_0(z)$ with  $\psi_0 \in A_1^2(\D)$, and the relation (8) is completed by
$$
\vk(f^{\mu_0}) = k(f^{\mu_0}) = q_L = \fc{1}{\rho_L} = \sup_{\psi \in A_1^2(D),\|\psi\|_{A_1} =1}
|\langle \mu_0, \psi\rangle_\D|,
$$
where $L = f^{\mu_0}(\mathbb S^1)$, and $q_L$ denotes its quasireflection coefficient.  }  (equal to the minimum of
dilatations $\partial_z W/\partial{- \ov z} W$ of the orientation reversing quasiconformal homeomorphisms of the sphere $\hC$ which preserve point-wise the quasicircle $L \subset \hC$ and interchange its interior and exterior domains.

\bigskip
This was established by the author earlier even for univalent functions in the arbitrary quasidisks (see, e.g.,
\cite{Kr6} (and in a special case by K\"{u}hnau \cite{Ku3}), but it follows from (7) immediately.

Recall that quasireflection coefficient of an (oriented) quasicircle $L^\prime \subset \hC$ is the minimumal dilatation  $\partial_z W/\partial{- \ov z} W$ of the orientation reversing quasiconformal homeomorphisms of the sphere $\hC$, which preserve point-wise the quasicircle $L \subset \hC$ and interchange its interior and exterior domains.

Theorem 1 also has other important consequences and applications. These results will be presented elsewhere.

\bigskip
\centerline{\bf 3. PROOF OF THEOREM 1}

\bigskip\noindent
$\mathbf {1^0}$.
First observe that the Grunsky coefficients $\a_{m n}(f^\mu)$ of functions $f^\mu \in \Sigma_Q$ generate for each $\x = (x_n) \in l^2$ with $\|\x\| = 1$ the holomorphic maps
 \be\label{9}
h_{\x}(f^\mu) = \sum\limits_{m,n=1}^\iy \ \a_{m n} (f^\mu) x_m x_n : \ \Belt(D)_1 \to \D,
\end{equation}
and
 \be\label{10}
\sup_{\x} |h_{\x}(f^\mu)| = \vk_{D^*}(f^\mu).
\end{equation}
The holomorphy of these functions follows from the holomorphy of coefficients $\a_{m n}$ with respect to
Beltrami coefficients $\mu \in \Belt(D)_1$ mentioned above, by applying the estimate
\be\label{11}
\Big\vert \sum\limits_{m=j}^M \sum\limits_{n=l}^N \ \beta_{mn} x_m x_n \Big\vert^2 \le \sum\limits_{m=j}^M |x_m|^2
\sum\limits_{n=l}^N |x_n|^2
\end{equation}
which holds for any finite $M, N$ and $1 \le j \le M, \ 1 \le l \le N$ (see \cite{Po}, p. 61).

Similar arguments imply that the maps (9) regarded as functions of points $\vp^\mu = S_{f^\mu}$ in the universal Teichm\"{u}ller space \ $\T$ are holomorphic on $\T$. This holomorphy provides, together
with the equality (10), that the Grunsky norm $\vk_{D^*}$ regarded as a function of the Schwarzians $S_f$ is
logarithmically plurisubharmonic on the space $\T$. In addition, as it is established in \cite{Kr6}, the functions
$\vk_{D^*}(S_f)$ is continuous; moreover, it is Lipschitz continuous on this space. The Teichm\"{u}ller norm
has the similar properties.

Now, given a function $f \in \Sigma_Q$, take its extremal extension $f^\mu$ (i.e., such that $k(f) = \|\mu\|_\iy$)
and consider its extremal disk
$$
\D(\mu) = \{t \mu/\|\mu\|_\iy: \ |t| < 1\} \subset \Belt(\D)_1.
$$
Put $\mu^* = \mu/\|\mu\|_\iy$.

Using the functions (9), we pull back the hyperbolic metric $\ld_\D(t)|dt| = |dt|/(1 - |t|^2)$ of the disk $\D$, which
implies on this disk the conformal metrics $\ld_{h_\x}(t) |dt|$ with
$$
\ld_{h_\x}(t) = |h_\x^\prime (t)|/(1 - |h_\x(t)|^2).
$$
All these metric have at their noncrical points the Gaussian curvature $-4$.
We take the upper envelope of these metrics
$$
\ld_\vk(t) = \sup \{\ld_{h_\x}(t): \ \x \in S(l^2)\}
$$
followed by its upper semicontinuous regularization, which determines a logarithmically
subharmonic metric $\ld_\vk(t)$ on the unit disk. It is circularly symmetric, i.e., $\ld_\vk(t) = \ld_\vk(|t|)$.
This metric being generated by the Grunsky coefficients of $f$ actually is the differential (infinitesimal) form
of the norm $\vk(f)$.

On a standard way, one obtains that $\ld_\vk$ has at any its noncritical point $t_0$  a supporting subharmonic metric $\ld_0$ of Gaussian curvature at most $- 4$, and hence, $\kappa_{\ld_\vk} \le - 4$; the details see, e.g., in \cite{Kr3}.

It is shown in \cite{Kr7}, that if the initial function $f \in \Sigma_Q$ has in $\D^*$ the expansion
$$
f(z) = z + b_n z^n + b_{n+1} z^{-n-1} + \dots, \quad b_n \ne 0, \ \ n \ge 1.
$$
then the corresponding metric $\ld_\vk$ satisfies the relation
  \be\label{12}
\lim\limits_{r\to 0} \fc{\ld_\vk(r)}{r^n}  = \sup_{\x \in S(l^2)} \fc{|d h_{\x}(0) (t \mu_0/\|\mu_0\|_\iy|)}{|t|} \ge  \sup_{\psi \in A_1^2(\D),\|\psi\|_{A_1} =1}  \ \Big\vert \iint\limits_\D \fc{\mu_0(z)}{\|\mu_0\|_\iy} \ \psi(z) dx dy\Big\vert,
\end{equation}
which is an infinitesimal version of (4).

\bigskip\noindent
$\mathbf {2^0}$.
To find an upper bound for $\vk(f)$(smaller that $k(f)$), we apply to functions $h_{\x}(f^{t\mu^*})$ a special case of Golusin's improvement \cite{Go} of the classical Schwarz lemma given by the following

\bigskip\noindent
{\bf Lemma 3}. {\it A holomorphic function
$$
g(t) = c_m t^m + c_{m+1} t^{m+1} + \dots: \D \to \D \quad (c_m \ne 0, \ \ m \ge 1),
$$
in $\D$ is estimated by
 \be\label{13}
|g(t)| \le |t|^m \fc{|t| + |c_m|}{1 + |c_m| |t|},
\end{equation}
with equality only for the function
 \be\label{14}
g_m(t) = t^m  \fc{t+ c_m}{1 + \ov c_m t} = c_m t^m + \dots \ .
\end{equation}
The corresponding differential metrics are related by
 \be\label{15}
g^* \ld_\D(r) \le \ld_j(r) := \fc{j^\prime(r)}{1 - j(r)^2}
\end{equation}
where $r = |t| < 1, \ld_\D(r) = 1/(1 - r^2)$ is the hyperbolic metric of the unit disk of curvature $-4$, and
$$
j(r) = r^m \ \fc{r + |c_m|}{1 + r |c_m|}.
$$
}

\bigskip
Note that the pull-backed metric $\ld_j(r)$ has the constant curvature $- 4$ in the punctured disk $\{0 < |t| < 1\}$
and in view of (10) any function (9) is dominated by (14) with $|c_m| \ge \a_\D(f)$.

\bigskip\noindent
$\mathbf {3^0}$.
To calculate explicitly the quantity $\a(f)$ given by the right hand side of (4), one can
use the variational formula for the functions
$f^\mu(z) = z + b_0 + b_1 z^{-1} + \dots$ from $\Sigma^0$ with extensions satisfying $f^\mu(0) = 0$ (cf \cite{Kr7}).
Namely, for small $\|\mu\|_\iy$,
$$
 f^\mu(z) = z - \fc{1}{\pi} \iint\limits_\D \mu(w) \left( \fc{1}{\zeta - z} -
\fc{1}{\zeta} \right) d\xi d\eta + O(\|\mu^2\|_\iy), \quad \zeta = \xi + i\eta,
$$
where the ratio $O(\|\mu^2\|_\iy^2)/\|\mu^2\|_\iy^2$ is uniformly bounded on compact sets of $\C$. Then
$$
b_n = \fc{1}{\pi} \iint\limits_\D \mu(\zeta) \zeta^{n-1} d\xi d\eta + O(\|\mu^2\|_\iy), \quad n = 1, 2, \dots,
$$
and
$$
\a_{m n}(\mu) = - \pi^{-1} \iint\limits_\D \mu(z) z^{m+n-2} dx dy + O(\|\mu\|_\iy^2), \quad
\|\mu\|_\iy \to 0.
$$
Using the normalized coefficient
$\mu^* = \mu/\|\mu\|_\iy$ and passing accordingly to coefficients $t \mu^*$ with $|t| < 1$, one derives that
the differential at zero of the corresponding map $h_\x(t \mu^*)$ with $\x = (x_n) \in S(l^2)$ along the extremal disk
$\{t \mu^* \in \Belt(\D)_1\}: \ |t| < 1\}$ (and its image in the space $\T$) is given by
$$
d h_{\x}(0) (t \mu^*) = - \fc{t}{\pi} \iint\limits_\D \mu^*(z) \sum\limits_{m+n=2}^\iy \sqrt{m n} \ x_m x_n z^{m+n-2}
dx dy
$$

On the other hand, as was established in \cite{Kr1}, the elements of $A_1^2(\D)$ are represented in the form
$$
\psi(z) = \om(z)^2 = \fc{1}{\pi} \sum\limits_{m+n=2}^\iy \ \sqrt{mn} \ x_m x_n z^{m+n-2} \quad \text{with} \ \
\|\x\|_{l^2} = \|\om\|_{L_2}.
$$

Now, letting $|t| < 1$ and applying the relation (4) and Lemma 3 to $h_\x(t \mu^*)$, one obtains
$$
|h_\x(t\mu^*)| \le |t| \fc{|t| + |\langle \mu^*, \psi\rangle_\D|}{1 + |\langle \mu^*, \psi\rangle_\D| |t|}
$$
and, passing to the supremum over $\x \in S(l^2)$,
 \be\label{16}
\vk(f^{t \mu^*}) \le |t| \fc{|t| + \a(\mu^*)}{1 + \a(\mu^*) |t|},
\end{equation}
with
$$
\a(\mu^*) = \sup_{\psi \in A_1^2(\D),\|\psi\|_{A_1} =1}  \ |\langle \mu^*, \psi\rangle_\D| =: \a.
$$
The relations (4), (12), (16) result in
 \be\label{17}
\a |t| \le \vk(f^{t \mu^*}) \le |t| \fc{|t| + \a}{1 + \a |t|}.
\end{equation}

\bigskip\noindent
$\mathbf {4^0}$. Now we may prove the theorem. We first establish the relation between the differential (infinitesimal)
metric $\ld_\vk$ corresponding to $\vk(f)$ and its dominant of type (15) defined in accordance to Lemma 3.

Noting that the right hand side fraction in (13) is increasing with respect to $|c_m|$ on $[0, 1]$, we take the functions
$$
g_a(t) = t^m (t+ a)/(1 + \ov a t)
$$
with $ a > \a$ and define by (15) the corresponding metrics
 \be\label{18}
\ld_a(r) = j_a^\prime(r)/(1 - j_a(r)^2)
\end{equation}
(all of curvature $- 4$). We compare these metrics with metric $\ld_\vk$ generated by the Grunsky coefficients
of $f^\mu$, as was described in the first step.

Comparison of $\ld_\vk$ with any metric (18) (which is maximal for a given value of $c_m = a$) provides by (12) the strict inequality
$$
\ld_\vk(r) < \ld_a(r), \quad r < 1,
$$
(for any $a > \a$). In the limit as $a \to \a$, one obtains
 \be\label{19}
\ld_\vk(r) \le \ld_\a(r) = \a  + (1 - \a^2) r + \dots \quad r < 1.
\end{equation}

On the other hand, as was establishes in \cite{Kr7}, we have at the origin $t = 0$ the equality
$$
\ld_\vk(0) = \wt \ld_\a(0),
$$
where
$$
\wt \ld_\a(r) = \a/(1 - \a^2 r^2) = \a + O(r^2), \quad r \to 0,
$$
is a supporting conformal metric for $\ld_\vk$ at the origin with constant Gaussian curvature $- 4$.
Therefore, for small $r$,
 \be\label{20}
\a \le \ld_\vk(r) \le \ld_\a(r) \le \a  + (1 - \a^2) r + O(r^2).
\end{equation}

Now choose a sufficiently small neighborhood $U_0$ of the origin $t = 0$ and put
$$
M = \{\sup \ld_\a(t): t \in U_0\}.
$$
Then in this neighborhood, we have $\ld_\a(t) + \ld_\vk(t) \le 2M$. Consider the function
$$
u = \log \fc{\ld_\vk}{\ld_\a}.
$$
Then (cf. \cite{Mi}, \cite{Di}) for $t \in U_0$,
$$
\Delta u(t) = \Delta \log \ld_\vk(t) - \Delta \log \ld_\a(t)
\ge 4 (\ld_\vk^2 - \ld_\a^2) \ge 8M (\ld_\vk - \ld_\a).  
$$
The elementary estimate
$$
M \log(t/s) \ge t - s \quad \text{for} \ \ 0 < s \le t < M
$$
(with equality only for $t = s$) implies that
$$
M \log \fc{\ld_\a(t)}{\ld_\vk(t)} \ge \ld_\a(t) - \ld_\vk(t),
$$
and hence,
$$
\Delta u(t) \ge 8 M^2 u(t).
$$
It follows that the function $u$ satisfies on $U_0$ the inequality (6) with $K = 8 M^2$.

Applying to $u$ Lemma 2 and noting that (20) yields
$$
u(0) = \lim\limits_{t \to 0} \log \fc{\ld_\vk(t)}{\ld_\a(t)} = 0,
$$
one derives by this lemma that both metrics $\ld_\vk$ and $\ld_\a$ must be equal on $U_0$, and in the similar way their equality on the entire disk $\D$.

\bigskip\noindent
$\mathbf {4^0}$.
Finally, we apply the following reconstruction lemma for the Grunsky norm proven in \cite{Kr3}, which provides that this norm is  the integrated form of $\ld_\vk$ along the Teichm\"{u}ller extremal disks.

\bigskip\noindent
{\bf Lemma 4}. {\it On any  Teichm\"{u}ller extremal disk
$\D(\mu_0) = \{t \mu_0/\|\mu_0\|_\iy: |t| < 1\}  \subset \Belt(\D)_1$,
we have the equality}
  \be\label{21}
\tanh^{-1}[\vk(f^{r\mu_0/\|\mu_0\|_\iy})] = \int\limits_0^r \ld_\vk(t) dt.
\end{equation}

\bigskip
Integrating the metrics  $\ld_\vk$ and $\ld_\a$ along the indicated extremal disks, one obtains from (21) and the right
equality in (4) the required equality (7), completing the proof of the theorem.

\bigskip
\centerline{\bf 4. EXAMPLE}

\bigskip
The simplest example illustrating Theorem 1 is given by the map
 \be\label{22}
f_{3,t}(z) = z \Bigl(1 + \fc{t}{z^3}\Bigr)^{2/3}, \quad 0 \le |t| < 1,
\end{equation}
whose extremal extension to $\D$ is
$$
\wh f_{3,t}(z) = z \Bigl[1 + t \Bigl(\frac{|z|}{z} \Bigr)^3\Bigr]^{2/3}
$$
with Beltrami coefficient
$\mu_3(z) := \mu_{\wh f_{3,t}}(z) = t|z|/z$.
In polar coordinates, $\mu_3(r e^{i \theta}) = t e^{- i \theta}$.

This map has threefold rotational symmetry
$$
f(e^{2n \pi i/3} z) = e^{2n \pi i/3} f(z), \quad n = 0, 1, 2.
$$
Note that $f_{3,t}(z) = f_t(z^{3/2})^{2/3}$, where
$$
f_t(z) = \begin{cases}  z + t/z, \ \ &|z| \ge 1,  \\
                        z + t \ov z, &|z| < 1.
\end{cases}
$$
For this function, we have (see \cite{Kr4})
 \be\label{23}
\sup_{\psi \in A_1^2, \|\psi\|_{A_1} = 1} \ |\langle \mu_3^*, \psi \rangle_\D|
= \max_{\psi_3} \ |\langle \mu_3^*, \psi_3 \rangle_\D|= \fc{2 \sqrt{2}}{3},
\end{equation}
where $\psi_3$ are the squares of nonconstant linear functions
$$
\psi_3(z) = \om_1(z)^2 := (a_0 + a_1 z)^2
$$
with $a_0 \ne 0, \ a_1 \ne 0$, and $\|\psi_3\|_{A_1} = 1$.

The map (22) was the first explicit example of functions $f \in \Sigma_Q$ with $\vk(f) < k(f)$. This was established
in 1981 independently by K\"{u}hnau in \cite{Ku2}, applying the technique of Fredholm eigenvalues, and by the author
in a different way (this also follows from the general result mentioned above).
The quantity (23) provides the value of reciprocal to eigenvalue $\rho_{f_{3,t}} (\mathbb S^1)$.

More generally, all functions $f_m(z) = f_t(z^{m/2})^{2/m}$ with odd $m = 2p - 1$ have similar features.

\bigskip
{\small\em{ \leftline{Department of Mathematics, Bar-Ilan
University, 5290002 Ramat-Gan, Israel} \leftline{and
Department of Mathematics, University of Virginia,  Charlottesville, VA 22904-4137, USA}}

\end{document}